\documentclass[reqno,12pt]{amsart}
\usepackage{epsf}\usepackage{amssymb} \usepackage{cite}

\newcommand{\beq}{\begin{equation}}
\newcommand{\ee}{\end{equation}}
\newcommand{\bea}{\begin{eqnarray}}
\newcommand{\eea}{\end{eqnarray}}

\usepackage{hyperref}
\def\stackreb#1#2{\ \mathrel{\mathop{#1}\limits_{#2}}}

\newcommand{\CC}{\mathbb C}
\newcommand{\R}{\mathbb R}
\newcommand{\Z}{\mathbb Z}

\begin{document}

\title[Rational hypergeometric identities]
{Rational hypergeometric identities}

\author{Gor A. Sarkissian$^{1,2,3}$ and Vyacheslav P. Spiridonov$^{1,2}$}

\address{
${}^1$ Laboratory of Theoretical Physics, JINR, Dubna, 141980,  Russia
\newline
\indent
${}^2$ St. Petersburg Department of the Steklov Mathematical Institute
of Russian Academy of Sciences, Fontanka 27, St. Petersburg, 191023 Russia
\newline
\indent
${}^3$ Department of Physics, Yerevan State University, Yerevan, Armenia
}

\begin{abstract}
A special singular limit $\omega_1/\omega_2\to1$ is considered for the Faddeev modular
quantum dilogarithm (hyperbolic gamma function) and the corresponding hyperbolic integrals.
It brings a new class of hypergeometric identities associated with bilateral
sums of Mellin-Barnes type integrals of particular Pochhammer symbol products.
\end{abstract}

\maketitle


\vspace{1.0em}

Let us consider Faddeev's modular quantum dilogarithm \cite{Fad94,Fad95} (also known as the
hyperbolic gamma function)
\beq
\gamma^{(2)}(u;\mathbf{\omega})= \gamma^{(2)}(u;\omega_1,\omega_2)=e^{-\frac{\pi{\rm i}}{2}
B_{2,2}(u;\mathbf{\omega}) } \gamma(u;\mathbf{\omega}),
\label{HGF}\ee
where $B_{2,2}(u;\mathbf{\omega})$ is the second-order multiple Bernoulli polynomial
$$
B_{2,2}(u;\mathbf{\omega})=\frac{1}{\omega_1\omega_2}
\left((u-\frac{\omega_1+\omega_2}{2})^2-\frac{\omega_1^2+\omega_2^2}{12}\right)
$$
and
\beq
\gamma(u;\mathbf{\omega})= \frac{(\tilde q e^{2\pi {\rm i} \frac{u}{\omega_1}};\tilde q)_\infty}
{(e^{2\pi {\rm i} \frac{u}{\omega_2}};q)_\infty}
=\exp\left(-\int_{\R+{\rm i} 0}\frac{e^{ux}}
{(1-e^{\omega_1 x})(1-e^{\omega_2 x})}\frac{dx}{x}\right);
\label{int_rep}\ee
here $(z;q)_\infty=\prod_{j=0}^\infty(1-zq^j)$,
$q= e^{2\pi {\rm i}\frac{\omega_1}{\omega_2}},$ and $\tilde q= e^{-2\pi {\rm i}\frac{\omega_2}{\omega_1}}$.
The infinite product representation in (\ref{int_rep}) is valid for $|q|<1$ (i.e., when ${\rm Im}(\omega_1/\omega_2)>0$) and defines an analytic continuation of $\gamma(u;\mathbf{\omega})$
to arbitrary $u\in\CC$.
The integral representation is well defined for a wider domain of values of $q$. Thus, if
${\rm Re}(\omega_1), {\rm Re}(\omega_2)>0$, then the integral in (\ref{int_rep}) converges
for $0<{\rm Re}(u)< {\rm Re}(\omega_1+\omega_2)$, so that for $\omega_1/\omega_2>0$, when $|q|=1$,
the function $\gamma^{(2)}(u;\mathbf{\omega})$ remains analytic. Note also that this function is
symmetric in $\omega_1$ and $\omega_2$, although this is not evident at all in
its infinite product form. For simplicity, we assume below that
arg$(\omega_1)>{\rm arg}(\omega_2)$ and ${\rm Re}(\omega_1), {\rm Re}(\omega_2)>0$.
Also, we fix $\sqrt{\omega_1\omega_2}>0$, which can always be done due to homogeneity:
$\gamma^{(2)}(\lambda u;\lambda\omega_1,\lambda\omega_2)=\gamma^{(2)}(u;\omega_1,\omega_2)$,
$\lambda\neq 0$.

It is not difficult to deduce the following asymptotic formulas:
\bea\nonumber
\makebox[-2em]{}
&& {\rm cone\; I}:\quad \stackreb{\lim}{z\to \infty}e^{{\pi{\rm i}\over 2}B_{2,2}(z;\omega_1,\omega_2)}\gamma^{(2)}(z;\omega_1,\omega_2)=1,
\quad {\rm arg}\;\omega_1<{\rm arg}\; z<{\rm arg}\;\omega_2+\pi,
\\ \nonumber \makebox[-2em]{}
&& {\rm cone\; II}: \quad\stackreb{\lim}{z\to \infty}e^{-{\pi{\rm i}\over 2}B_{2,2}(z;\omega_1,\omega_2)}\gamma^{(2)}(z;\omega_1,\omega_2)=1,
\quad  {\rm arg}\;\omega_1-\pi<{\rm arg}\; z<{\rm arg}\;\omega_2.
\eea

Recently, in \cite{GSVS}, we found a new singular degeneration of the function
$\gamma^{(2)}(u;\mathbf{\omega})$ emerging in the limit
$\sqrt{\omega_1/\omega_2}=1+{\rm i}\delta$, $\delta\to 0^+$.
For generic values of the argument $u$, there are no singularities, but
under a special choice of $u$, uniformly on compact sets, one has
\beq\label{gam2b1}
\gamma^{(2)}\Big(\sqrt{\omega_1\omega_2}(n+1+y\delta);\mathbf{\omega}\Big)
\stackreb{=}{\delta\to 0^+} e^{-\frac{\pi {\rm i}}{2}n^2}
(4\pi\delta)^n\Big(\frac{1-n-{\rm i}y}{2}\Big)_n,
\ee
where $n\in \Z,\, y\in\CC,$ and $(a)_n$ is the standard Pochhammer symbol
$$
(a)_n=\frac{\Gamma(a+n)}{\Gamma(a)}=\left\{
\begin{array}{cl}
 a(a+1)\cdots(a+n-1), \quad & {\rm for} \; n>0,\quad \\
{\textstyle 1 \over \textstyle (a-1)(a-2)\cdots(a+n)}, \quad &{\rm for} \; n<0,
\end{array}
\right.
$$
So, $\gamma^{(2)}(u;\mathbf{\omega})$ tends to vanish for $n>0$ and to blow up for $n<0$.
This is somewhat similar to the limit $\omega_1+\omega_2\to 0$, also discussed
in \cite{GSVS}. In this paper we apply the limit (\ref{gam2b1})
to some univariate hyperbolic integral identities.

The elliptic beta integral \cite{spi:umn} is the most general univariate integral of hypergeometric type
admitting exact evaluation. As shown in \cite{rai:limits}, it can be rigorously
degenerated to the computable hyperbolic integral
\begin{eqnarray}
&& \int_{-{\rm i}\infty}^{{\rm i}\infty}
\Delta(z)\frac{dz}{2{\rm i}\sqrt{\omega_1\omega_2}}
= \prod_{1\leq j<k\leq 6}\gamma^{(2)}(g_j+g_k;\mathbf{\omega}),
\label{hyper}\end{eqnarray}
where
\beq
\Delta(z)=\frac{\prod_{k=1}^6\gamma^{(2)}(g_k+ z;\mathbf{\omega})\gamma^{(2)}(g_k- z;\mathbf{\omega}) }
{\gamma^{(2)}(2 z;\mathbf{\omega})\gamma^{(2)}(- 2 z;\mathbf{\omega})}
\label{kernel}\ee
and the balancing condition $\sum_{k=1}^6 g_k=\omega_1+\omega_2$ holds.

The poles and zeros of the function $\gamma^{(2)}(u;\mathbf{\omega})$ are located at
$$
u_{poles}\in \{ -n\omega_1 -m\omega_2\},\quad u_{zeros}= (n+1)\omega_1+(m+1)\omega_2,
\quad  n,m \in \Z_{\geq0},
$$
so that $\Delta(z)$ has the poles
\beq
z_{poles}\in \{ g_k+n\omega_1+ m\omega_2\}\cup \{ -g_k-n\omega_1- m\omega_2,\},\quad n,m \in \Z_{\geq0},
\quad k=1,\ldots, 6.
\label{IP}\ee
The $z\to\infty$ asymptotics of $\Delta(z)$  has the form
\beq
{\rm cone\; I}: \quad e^{6\pi{\rm i}\frac{z(\omega_1+\omega_2)}{\omega_1\omega_2}}, \qquad
{\rm cone\; II}: \quad e^{-6\pi{\rm i}\frac{z(\omega_1+\omega_2)}{\omega_1\omega_2}},
\label{AB}\ee
i.e., $\Delta(z)$ vanishes exponentially fast  for $z$ in these cones.
Formula (\ref{hyper}) is true for any integration contour separating the two arrays of poles in (\ref{IP})
and asymptotically lying inside the specified cones.

Now we set $\sqrt{\omega_1/\omega_2}=1+{\rm i}\delta$, fix
\beq
g_k=\sqrt{\omega_1\omega_2}(N_k+{\rm i}a_k\delta), \quad N_k\in \Z+\nu,\quad
\nu=0, \frac{1}{2}, \; a_k\in\CC,
\label{g}\ee
and consider the limit as $\delta\to 0^+$ (the origin of the parameter $\nu$ will be explained later).
One can verify that, for any choice of the discrete variables $N_k$, a number of poles start to collide
in pairs, pinching any admissible integration contour from two sides.
For Re$(g_k)>0$, we took the imaginary axis as the integration contour, but even that is
inevitably pinched by pairs of poles. One can estimate the rate of divergence associated with this
pinching by computing residues, but we choose instead a different strategy.
We evade pinching by scaling the integration variable, $z \propto y\delta$, so that there remains nontrivial
integration over $y$, although it takes a much simpler form.

Note that, for $\delta\to 0^+$, the balancing condition $\sum_{k=1}^6 g_k=\omega_1+\omega_2$ reduces to
two separate constraints on the discrete and continuous variables:
\beq\label{xknk2}
\sum_{k=1}^6 N_k=2, \qquad \sum_{k=1}^6 a_k=0.
\ee
If at least one $N_k$ is negative, it is necessary to deform the contour of integration to a Mellin-Barnes
type one in order to separate the corresponding sequences of poles.

The integral kernel $\Delta(z)$ starts to blow up when $z/\sqrt{\omega_1\omega_2}$ approaches
an integer and half-integer values, i.e., for
\beq
z=\sqrt{\omega_1\omega_2}(N+{\rm i}y\delta), \quad N\in\Z+\nu, \quad \nu=0, \frac{1}{2}, \; y\in\CC.
\label{z}\ee
Applying relation (\ref{gam2b1}), we obtain
$$
\Delta(z)\stackreb{\to}{\delta\to 0^+} (-1)^{2\nu} \frac{y^2-N^2}{(4\pi\delta)^6}
\prod_{k=1}^6 \Big(1+\frac{a_k-N_k\pm(y-N)}{2}\Big)_{N_k-1\pm N},
$$
where we use the convention $(a\pm b)_{n\pm m}=(a+b)_{n+m}(a-b)_{n-m}$.
The new discrete parameter $\nu=0, \frac{1}{2}$ in this formula emerged from the
fact that only the combinations $N_k\pm N$ should be integers, which evidently produces
the indicated freedom in $\nu$-values.

As to the right-hand side expression in (\ref{hyper}), it has the asymptotics
$$
RHS\stackreb{\to}{\delta\to 0^+}\frac{(-1)^{2\nu} {\rm i}}{(4\pi\delta)^5}
\prod_{1\leq j< k \leq 6}\Big(1+\frac{a_j+ a_k-N_j-N_k}{2}\Big)_{N_j+N_k-1}.
$$
As a result of these manipulations, in the limit as $\delta\to0^+$, relation (\ref{hyper})
reduces to the following exact formula.

{\bf Theorem 1.}
For parameters $N_k\in\Z$ and $a_k\in\CC$ satisfying constraints (\ref{xknk2}), one has
 \bea\nonumber
\frac{1}{8\pi{\rm i}}\sum_{N\in \mathbb{Z} +\nu}\int_{C_N} (y^2-N^2)
\prod_{k=1}^6 \Big(1+\frac{a_k-N_k\pm(y-N)}{2}\Big)_{N_k-1\pm N} dy
\\ \makebox[4em]{}
=\prod_{1\leq j< k \leq 6}\Big(1+\frac{a_j+ a_k-N_j-N_k}{2}\Big)_{N_j+N_k-1},
\label{ratbeta}\eea
where the contour $C_N$ separates the poles emerging from the Pochhammer symbols
corresponding to the plus sign (lying to the left of the contour)
from the poles corresponding to the factors  with the minus sign (lying to the right of the contour).

\smallskip

Formally, we wrote the infinite bilateral sum over $N$ of terms corresponding to all admissible cases, without specifying the exact values of this discrete variable and the contours $C_N$ that indeed make a
contribution to the sum. One can see that, for sufficiently
large $|N|$, all poles of the integrand lie on one side of the integration contour
and thus all of the corresponding integrals vanish, i.e., the sum over $N$ is finite.

The asymptotics of the Pochhammer symbol
$$
(z)_a={\Gamma(z+a)\over \Gamma(z)}\stackreb{=}{z\to\infty} z^a(1+o(1)), \quad a\in \CC,
$$
shows that all integrals in (\ref{ratbeta}) converge, because
the absolute values of the integrands decay when $|N+{\rm i}y|\to \infty$ as $|N+{\rm i}y|^{-6}$
(i.e., even the formal infinite bilateral sum of integrals converges).

Now consider the three simplest cases illustrating the nature of the obtained rational hypergeometric identities.
Suppose that $N_j=0$ and Im$(a_j)<0$ for $j=1,\ldots, 4$, $N_{5,6}=1$ and Im$(a_5+a_6)>0$.
For $N>0$ or $N<0$, all poles of the integrand lie on one side of the integration contour,
and therefore only the term $N=0$ gives nontrivial contribution to the sum over $N$.
The parameters $a_{5,6}$ drop out completely from the $N=0$ integral,
and the resulting formula takes the form (after dropping the common factor $2^5$)
$$
\frac{1}{\pi{\rm i}}\int_{\mathbb{R}}\frac{y^2dy}{\prod_{j=1}^4(a_j\pm y)} =
\sum_{j=1}^4\frac{a_j}{\prod_{k=1,\neq j}^4(a_k^2-a_j^2)}
=-\frac{a_1+a_2+a_3+a_4}{\prod_{1\leq j<k\leq 4}(a_j+a_k) }.
$$
Here the original integration contour (the imaginary axis) was turned into the real $y$-axis
and all $a_j$ lie below that axis. Closing the contour in the upper half-plane, we compute
the integral as the sum of the four pole residues $y=-a_j$. The last equality can be checked by a direct
computation of the sum. We remark that the derived formula shows the evaluation
of the leading $a_{j}\to 0$ asymptotics of the de Branges-Wilson beta integral \cite{aar}.

Suppose that $N_j=0$ and Im$(a_j)<0$ for $j=1,\ldots,5,$ $N_6=2,$ and Im$(a_6)>0$. Again only the term with $N=0$
gives a contribution, which results in the relation
$$
\frac{1}{\pi{\rm i}}\int_{\mathbb{R}}\frac{y^2(a_6\pm y)dy}{\prod_{j=1}^5(a_j\pm y)} =
\sum_{j=1}^5\frac{a_j(A^2-a_j^2)}{\prod_{k=1,\neq j}^5(a_k^2-a_j^2)}
=\frac{\prod_{j=1}^5(a_j-A)}{\prod_{1\leq j<k\leq 5}(a_j+a_k)}, \quad
A=\sum_{j=1}^5a_j.
$$
This formula corresponds to the $a_j\to 0$ asymptotics of the Rahman beta integral \cite{aar}.

Suppose that $N_j=\frac{1}{2}$ for $j=1,\ldots, 5$ and $N_6=-\frac{1}{2}$, and that the $a_j$ take arbitrary generic values with $\sum_{j=1}^6a_j=0$.
Nonzero terms in the sum emerge only for the two values $N=\pm \frac{1}{2}$, which give equal contributions.
As a result, we come to the identity
$$
\frac{2}{\pi{\rm i}}\int_{\mathbb{R}}\frac{(y^2-\frac{1}{4})dy}{(a_6-y)((a_6+y)^2-1)\prod_{j=1}^5(a_j+y)}
=\frac{1}{\prod_{j=1}^5(a_j+a_6)},\quad a_6=-\sum_{k=1}^5a_k,
$$
where only the pole $y=a_6$ lies below $\mathbb{R}$. An independent proof of this identity is very
simple---it reduces to the computation of the $y=a_6$ pole residue.
Perhaps, this and all other cases of identity (\ref{ratbeta}) correspond to
some yet-to-be-found transcendental plain hypergeometric beta integrals.

Consider now the hyperbolic analogue of the Euler-Gauss hypergeometric function, namely the integral
\beq
I(\underline{g})=\int_{-{\rm i}\infty}^{{\rm i}\infty}{\prod_{j=1}^8\gamma^{(2)}(g_j\pm z;\omega_1,\omega_2)\over
\gamma^{(2)}(\pm 2z;\omega_1,\omega_2)}\frac{dz}{2{\rm i}\sqrt{\omega_1\omega_2}}
\ee
with the parameters $g_j$ satisfying the conditions ${\rm Re}(g_j)>0$ and
$\sum_{j=1}^8 g_j=2(\omega_1+\omega_2).$  As shown
in \cite{rai:limits}, $I(\underline{g})$ represents a hyperbolic degeneration of the
elliptic hypergeometric $V$-function \cite{spi:essays}. Thus, it inherits the
corresponding $W(E_7)$-group of symmetry transformations generated by the relation
\beq\label{ide1}
I(\underline{g})=\prod_{1\leq j< k \leq 4}\gamma^{(2)}(g_j+g_k;\omega_1,\omega_2)
\prod_{5\leq j< k \leq 8}\gamma^{(2)}(g_j+g_k;\omega_1,\omega_2)\, I(\underline{\lambda}),
\ee
where
\beq
\lambda_j=g_j+\xi, \quad \lambda_{j+4}=g_{j+4}-\xi, \quad j=1,2,3,4,\quad
\xi={1\over 2}(\omega_1+\omega_2-\sum_{j=1}^4 g_j).
\ee

Parametrization (\ref{g}) in the limit $\delta\to 0^+$ yields the balancing condition
\begin{equation}
\sum_{k=1}^8 N_k=4, \qquad \sum_{k=1}^8 a_k=0.
\label{balance3}\end{equation}
We set $X=\sum_{j=1}^4 a_j, \, L=\sum_{j=1}^4 N_j$, and, for $j=1,\ldots,4$,
\beq
M_j=N_j+1-\frac{L}{2},\, M_{j+4}=N_{j+4}-1+\frac{L}{2},\, s_j=a_j-\frac{X}{2},\,
s_{j+4}=a_{j+4}+\frac{X}{2}.
\label{M}\ee
Then
$$
{\xi\over \sqrt{\omega_1\omega_2}}=1-{L\over 2}-{\rm i}\delta{X\over 2},\quad
 {\lambda_k\over \sqrt{\omega_1\omega_2}}=M_k+{\rm i}s_k\delta,\quad k=1,\ldots, 8.
$$
The passage to the $\delta\to 0^+$ limit in relation (\ref{ide1}) yields the following statement.

{\bf Theorem 2.} For $N_k\in\Z $ and $a_k\in\CC$ satisfying constraints (\ref{balance3})
and $M_k$ and $s_k$, $k=1,\ldots, 8$, defined in (\ref{M}), one has  the identity
\bea \nonumber &&
\sum_{N\in\Z+\nu}\int_{C_N}(y^2-N^2)
\prod_{k=1}^8 \Big(1+\frac{a_k-N_k\pm(y-N)}{2}\Big)_{N_k-1\pm N} dy
\\  \nonumber && \makebox[3em]{}
=(-1)^L\prod_{1\leq j< k \leq 4\atop 5\leq j< k \leq 8}\Big(1+\frac{a_j+ a_k-N_j-N_k}{2}\Big)_{N_j+N_k-1}
\\   && \makebox[4em]{} \times
\sum_{M\in\Z+\mu}\int_{C_M}(y^2-M^2)
\prod_{k=1}^8 \Big(1+\frac{s_k-M_k\pm(y-M)}{2}\Big)_{M_k-1\pm M} dy,
\label{rat-trafo}\eea
where the contours $C_N$ and $C_M$ are defined as in Theorem 1 and $\mu, \nu=0, \frac{1}{2}$ are
different discrete parameters whose values depend on the parity of the integer $L$. If $L$ is even, then $\mu=\nu$, and if $L$ is odd, then $\mu \neq \nu$.

\smallskip

The symmetry transformation (\ref{rat-trafo}) extends the manifest $S_8$ permutation symmetry of
the parameters to $W(E_7)$---the Weyl group of the exceptional root system $E_7$.
By the same reason as in Theorem 1, the sums over $N$ and $M$
in relation (17) contain a finite number of terms. Therefore, on both sides of (\ref{rat-trafo}),
we have rational functions of $a_k$, which can be explicitly computed by methods of residue calculus;
i.e., we have a realization of the $W(E_7)$ group action on rational functions.

There are very many multiple hyperbolic integrals either obeying symmetry transformations or admitting exact
evaluation, which can be derived from limits of elliptic hypergeometric integrals \cite{rai:limits,spi:essays}.
All of them can be reduced as described above and yield multiple rational hypergeometric
identities. As to applications of the derived exact formulas, it would be interesting to give them
the form of star-triangle or star-star relations and connect to the Yang-Baxter equation.
Also they should have applications to two-dimensional vortex dynamics and two-dimensional
Liouville field theory with central charge $c=25$.

\medskip

{\bf Acknowledgments.} This work is supported by the Russian Science Foundation (project no. 19-11-00131).

\end{document}